\documentclass{bourbaki}
\usepackage[francais]{babel}

\date{Juin 2004}
\bbkannee{56{\`e}me ann{\'e}e, 2003-2004}
\bbknumero{937}
\title{{\'E}tats quasi-libres libres et facteurs de type III}
\subtitle{[d'apr{\`e}s D. Shlyakhtenko]}

\author{Stefaan VAES}
\address{CNRS, Institut de Math{\'e}matiques de Jussieu \\
Alg{\`e}bres d'Op{\'e}rateurs, Plateau 7E \\
175 rue du Chevaleret \\
75013 Paris}
\email{vaes@math.jussieu.fr}

\newcommand{\R}{\mathbb{R}}
\newcommand{\F}{\mathbb{F}}
\newcommand{\C}{\mathbb{C}}
\newcommand{\N}{\mathbb{N}}
\newcommand{\ot}{\otimes}
\newcommand{\recht}{\rightarrow}

\newcommand{\cst}{C$^*$}
\newcommand{\om}{\omega}
\newcommand{\si}{\sigma}
\newcommand{\Om}{\Omega}
\newcommand{\cF}{\mathcal{F}}
\newcommand{\vfi}{\varphi}
\newcommand{\free}[1]{\mbox{$\Gamma(#1)^{\prime\prime}$}}
\newcommand{\dpr}{^{\prime\prime}}
\newcommand{\freep}[1]{\underset{#1}{\mbox{\large $\ast$}}}
\newcommand{\id}{\text{id}}
\newcommand{\la}{\lambda}
\newcommand{\B}{\operatorname{B}}
\newcommand{\eps}{\varepsilon}
\newcommand{\cC}{\mathcal{C}}
\newcommand{\cE}{\mathcal{E}}
\newcommand{\al}{\alpha}
\newcommand{\zetatil}{\tilde{\zeta}}
\newcommand{\etatil}{\tilde{\eta}}
\newcommand{\gatil}{\tilde{\gamma}}

\newcommand{\Hbol}{{\overset{\circ}{H}}}
\newcommand{\Honebol}{{\overset{\circ}{H}_1}}
\newcommand{\Htwobol}{{\overset{\circ}{H}_2}}
\newcommand{\be}{\beta}
\newcommand{\ga}{\gamma}
\newcommand{\xbol}{{\overset{\circ}{x}}}
\newcommand{\etapr}{\eta^{\prime}}
\newcommand{\etadpr}{\eta^{\prime\prime}}
\newcommand{\mupr}{\mu^{\prime}}
\newcommand{\mudpr}{\mu^{\prime\prime}}
\newcommand{\zetapr}{\zeta^{\prime}}
\newcommand{\Aut}{\operatorname{Aut}}
\newcommand{\Out}{\operatorname{Out}}
\newcommand{\Ad}{\operatorname{Ad}}
\newcommand{\Int}{\operatorname{Int}}

\newcommand{\T}{\mathbb{T}}
\newcommand{\cS}{\mathcal{S}}

\newcommand{\car}{\operatorname{CAR}}

\newcommand{\libre}{\Gamma}
\newcommand{\cR}{\mathcal{R}}
\newcommand{\te}{\mathcal{\theta}}
\newcommand{\Z}{\mathbb{Z}}
\newcommand{\trace}{\operatorname{tr}}
\newcommand{\cliff}{\operatorname{Cliff}}

\renewcommand{\lq}{\og}
\renewcommand{\rq}{\fg}

{\theoremstyle{definition}\newtheorem{nota}[defi]{Notation}
\newtheorem{term}[defi]{Terminologie}
}

\begin{document}
\maketitle

Le but de cet expos{\'e} est de pr{\'e}senter une famille d'alg{\`e}bres de von Neumann introduite par Shlyakhtenko et de donner un aper{\c c}u des r{\'e}sultats de
classification des alg{\`e}bres de cette famille. Ces alg{\`e}bres de von Neumann sont construites dans le cadre des probabilit{\'e}s libres de Voiculescu.

Murray et von Neumann ont initi{\'e} la classification des alg{\`e}bres de von Neumann. Ils ont d{\'e}montr{\'e} que chaque alg{\`e}bre de von Neumann s'{\'e}crit comme
int{\'e}grale directe de facteurs (alg{\`e}bres de von Neumann de centre trivial) et ils ont classifi{\'e} les facteurs en diff{\'e}rents types~: I, II et III. Dans sa
th{\`e}se \cite{Connes-class}, Connes a raffin{\'e} cette classification en introduisant les sous-types III$_\la$ ($0 \leq \la \leq 1$). La construction de
Shlyakhtenko donne tout un monde d'exemples de facteurs de type
III$_1$, c'est-{\`a}-dire le plus haut dans la \lq hi{\'e}rarchie\rq\ des facteurs.

L'id{\'e}e de Shlyakhtenko est de donner dans le cadre des probabilit{\'e}s libres de Voiculescu
une version du foncteur $\car$ (relations d'anticommutation canoniques) et des {\'e}tats quasi-libres associ{\'e}s. Le foncteur $\car$ associe {\`a} chaque espace de Hilbert $H$ la \cst-alg{\`e}bre unif{\`e}re universelle
$\car(H)$ engendr{\'e}e par la famille $\{a(\xi) \mid \xi \in H\}$ telle que
\begin{enumerate}
\item $\xi \mapsto a(\xi)$ est lin{\'e}aire,
\item les relations d'anticommutation canoniques sont v{\'e}rifi{\'e}es
\begin{align*}
a(\eta) a(\xi)^* + a(\xi)^* a(\eta)  &= \langle \eta,\xi \rangle 1 \; ,\\
a(\xi) a(\eta) + a(\eta) a(\xi) &= 0 \; .
\end{align*}
\end{enumerate}
Les relations d'anticommutation canoniques peuvent {\^e}tre r{\'e}alis{\'e}es par les op{\'e}rateurs de cr{\'e}ation sur \emph{l'espace de Fock antisym{\'e}trique}. Tout
op{\'e}rateur $S$ agissant sur $H$ tel que $0 \leq S \leq 1$ donne lieu {\`a} un {\'e}tat $\om_S$ de la \cst-alg{\`e}bre $\car(H)$, qu'on appelle \emph{{\'e}tat
quasi-libre de covariance $S$}. Les repr{\'e}sentations GNS \cite{dixmier} associ{\'e}es aux {\'e}tats quasi-libres ont {\'e}t{\'e} beaucoup {\'e}tudi{\'e}es
\cite{shale-stine,powers-stormer,araki}. Dans une telle repr{\'e}sentation, le bicommutant $\car(H)\dpr$ est une alg{\`e}bre de von Neumann. Chaque {\'e}tat
quasi-libre sur $\car(H)$ donne donc une alg{\`e}bre de von Neumann qui s'av{\`e}re {\^e}tre un facteur. Les travaux d'Araki \& Woods \cite{araki-woods} et de
Powers \& St\o rmer \cite{powers-stormer} permettent de d{\'e}terminer le type de ces \emph{facteurs d'Araki-Woods}\footnote{Il d{\'e}coule des
  travaux de Powers \& St{\o}rmer que les facteurs associ{\'e}s aux {\'e}tats
  quasi-libres sont des produits tensoriels infinis de matrices $2$
  fois $2$ (des facteurs ITPFI$_2$). Plus g{\'e}n{\'e}ralement, Araki \& Woods
  {\'e}tudient et d{\'e}terminent le type des produits tensoriels infinis de facteurs de type I (les
  ITPFI) et il y a des ITPFI qui ne sont pas ITPFI$_2$.}. Plus pr{\'e}cis{\'e}ment, si $0 < \la \leq 1$, il existe
exactement une classe d'isomorphisme de facteurs d'Araki-Woods de type III$_\la$ et il existe une famille non-d{\'e}nombrable de facteurs d'Araki-Woods de
type III$_0$ et mutuellement non-isomorphes.

Remarquons que les facteurs d'Araki-Woods sont des
\emph{facteurs moyennables}. Dans \cite{Connes-class-injective} Connes
a r{\'e}ussi {\`a} donner une classification
compl{\`e}te des facteurs moyennables. C'est un des r{\'e}sultats les plus profonds
de la th{\'e}orie des alg{\`e}bres de von Neumann. Pour chacun des types
II$_1$, II$_\infty$, III$_\la$ ($0 < \la \leq 1$) il existe un unique
facteur moyennable. Les facteurs moyennables de type III$_0$ sont
classifi{\'e}s par un invariant en th{\'e}orie ergodique qu'on appelle \emph{flot
des poids}. L'unicit{\'e} du facteur moyennable de type III$_1$
est d{\^u} {\`a} Haagerup \cite{Haag-class}.
Dans le cadre des probabilit{\'e}s libres, nous allons obtenir des
facteurs non-moyennables. En particulier Shlyakhtenko construit une famille
non-d{\'e}nombrable de facteurs de type III$_1$, non-moyennables et
mutuellement non-isomorphes.

L'analogue en probabilit{\'e}s libres du foncteur $\car$ associe {\`a} un
espace de Hilbert r{\'e}el $H_\R$, la \cst-alg{\`e}bre universelle
$\libre(H_\R)$ engendr{\'e}e
par la famille $\{ s(\xi) \mid \xi \in H_\R \}$ telle que
\begin{itemize}
\item $s(\xi)$ est auto-adjoint pour tout $\xi \in H_\R$,
\item $\xi \mapsto s(\xi)$ est $\R$-lin{\'e}aire,
\item $\|s(\xi)\| \leq \|\xi\|$ pour tout $\xi \in H_\R$.
\end{itemize}
Pour chaque plongement isom{\'e}trique $H_\R \hookrightarrow H$ de $H_\R$
dans un espace de Hilbert $H$, on obtient une repr{\'e}sentation de
$\libre(H_\R)$ sur \emph{l'espace de Fock plein}
$$\cF(H) = \C \Om \oplus \bigoplus_{n=1}^\infty H^{\ot n} \; ,$$
en posant $s(\xi) =
(\ell(\xi) + \ell(\xi)^*)/2$ o{\`u} $\ell(\xi)$ est l'op{\'e}rateur de
cr{\'e}ation.

La construction des facteurs d'Araki-Woods libres suppose la donn{\'e}e d'un groupe {\`a} un param{\`e}tre $(U_t)$ de transformations orthogonales d'un espace de
Hilbert r{\'e}el $H_\R$ qui induit un plongement isom{\'e}trique $H_\R \hookrightarrow H$ de $H_\R$ dans le complexifi{\'e} $H$. On obtient une repr{\'e}sentation de
$\libre(H_\R)$ dont l'image est not{\'e} $\libre(H_\R,U_t)$. Le facteur d'Araki-Woods libre $\free{H_\R,U_t}$ est l'alg{\`e}bre de von Neumann engendr{\'e}e par
$\libre(H_\R,U_t)$. La restriction de l'{\'e}tat du vide au facteur $\free{H_\R,U_t}$ est \emph{l'{\'e}tat quasi-libre libre} not{\'e} $\vfi_U$. Alors,
$\free{H_\R,U_t}$ est un facteur de type III, sauf si $U_t = \id$ pour tout $t \in \R$.

Nous commen{\c c}ons cet expos{\'e} par rappeler la classification des facteurs de Connes et les probabilit{\'e}s libres de Voiculescu et par une
pr{\'e}sentation de plusieurs points de vue sur nos donn{\'e}es essentielles, les repr{\'e}sentations orthogonales de $\R$. Au \S\ref{sec.arakiwoodslibre} nous
d{\'e}finissons l'alg{\`e}bre de von Neumann $\free{H_\R,U_t}$ avec l'{\'e}tat quasi-libre libre $\vfi_U$. Nous {\'e}tudions en d{\'e}tail le cas $H_\R = \R^2$ muni de
la repr{\'e}sentation de $\R$ par rotations. Au \S\ref{sec.classification} nous pr{\'e}sentons les principaux r{\'e}sultats de classification et de
non-isomorphisme des facteurs d'Araki-Woods libres obtenus par Shlyakhtenko~:
\begin{itemize}
\item Classification compl{\`e}te des facteurs d'Araki-Woods libres associ{\'e}s {\`a} une repr{\'e}\-sen\-ta\-tion $(U_t)$ \emph{presque-p{\'e}riodique}.
\item Construction d'une famille non-d{\'e}nombrable de facteurs d'Araki-Woods libres non-presque-p{\'e}riodiques et mutuellement non-isomorphes. Ces facteurs
sont distingu{\'e}s par leur invariant $\tau$ de Connes.
\item Construction de deux facteurs d'Araki-Woods libres non-isomorphes
  et ayant le m{\^e}me invariant $\tau$ de Connes. Ceci est une application
de la notion de \emph{dimension entropique libre} introduite par Voiculescu \cite{Voic-entropy3,Voic-entropy2}.
\item D{\'e}monstration que la classe d'isomorphisme d'un facteur
  Araki-Woods libre peut d{\'e}pendre de la multiplicit{\'e} de la repr{\'e}sentation $(U_t)$. Ce r{\'e}sultat
est d{\'e}montr{\'e} {\`a} l'aide de la notion \emph{d'alg{\`e}bre de von Neumann solide} due {\`a} Ozawa \cite{ozawa}.
\end{itemize}

La classification compl{\`e}te des facteurs d'Araki-Woods libres reste un probl{\`e}me ouvert. Les r{\'e}sultats de non-isomorphisme pr{\'e}sent{\'e}s au
\S\ref{sec.classification} montrent qu'un facteur d'Araki-Woods libre $\free{H_\R,U_t}$ d{\'e}pend fortement de la classe de la mesure spectrale de la
repr{\'e}sen\-tation $(U_t)$. Ce probl{\`e}me de classification est beaucoup plus difficile que la classification des facteurs d'Araki-Woods, pour la raison
suivante. Powers et St\o rmer \cite{powers-stormer} d{\'e}montrent essentiellement que deux facteurs d'Araki-Woods associ{\'e}s {\`a} des {\'e}tats quasi-libres sont
isomorphes si leurs op{\'e}rateurs de covariance diff{\`e}rent d'un op{\'e}rateur d'Hilbert-Schmidt. En particulier, d'apr{\`e}s un r{\'e}sultat de von Neumann, il
suffit de consid{\'e}rer le cas d'un op{\'e}rateur de covariance diagonalisable. Ceci n'est plus le cas pour les facteurs d'Araki-Woods libres. Une grande
partie des facteurs d'Araki-Woods libres ne peut {\^e}tre obtenue par des repr{\'e}sentations orthogonales presque-p{\'e}riodiques.

Dans le dernier \S\ref{sec.appendice} nous pr{\'e}sentons un nouveau r{\'e}sultat sur les produits libres ce qui permet au \S\ref{subsec.types} de d{\'e}montrer
un r{\'e}sultat un peu plus g{\'e}n{\'e}ral que dans l'article \cite{shl3}.

Il y a un certain nombre de r{\'e}sultats et d'applications dans la
th{\'e}orie des facteurs d'Araki-Woods libres dont on ne parlera pas en
d{\'e}tail dans cet expos{\'e}. Notons que Shlyakhtenko a d{\'e}montr{\'e} dans
\cite{shl4} que les facteurs d'Araki-Woods libres $T_\la$ de type
III$_\la$ ($0 < \la < 1$) sont des \emph{facteurs premiers}~: ils ne peuvent
{\^e}tre {\'e}crits comme produit tensoriel de deux facteurs diffus (sans
projecteurs minimaux). Dans \cite{pisier-shl} Pisier et Shlyakhtenko
utilisent des facteurs d'Araki-Woods libres comme mod{\`e}les pour d{\'e}montrer
une \emph{in{\'e}galit{\'e} de Grothendieck} pour les espaces
d'op{\'e}rateurs. Dans \cite{Vouter} les facteurs d'Araki-Woods libres sont
utilis{\'e}s pour construire des \emph{actions ext{\'e}rieures} de groupes
quantiques localement compacts.

{\it Je remercie S.~Baaj, E.~Germain, D.~Shlyakhtenko et G.~Skandalis
  pour leur aide pendant la pr{\'e}paration de cet expos{\'e}.}

\section{Rappels} \label{sec.rappels}

\subsection{Alg{\`e}bres de von Neumann. Type des facteurs.}

On rappelle que pour un ensemble $X$ d'op{\'e}rateurs born{\'e}s sur un espace de Hilbert $H$, $X \subset \B(H)$, on appelle \emph{commutant} de $X$ et on
note $X'$ l'ensemble de tous les op{\'e}rateurs $T \in \B(H)$ qui commutent {\`a} $X$. Si $A \subset \B(H)$ est une sous-alg{\`e}bre involutive qui agit d'une
fa{\c c}on non-d{\'e}g{\'e}n{\'e}r{\'e}e sur $H$, le \emph{bicommutant} $A\dpr$ co{\"\i}ncide avec l'adh{\'e}rence de $A$ dans $\B(H)$ pour la topologie faible, c'est-{\`a}-dire la
topologie donn{\'e}e par les semi-normes $T \mapsto |\langle T \xi,\eta \rangle|$, o{\`u} $\xi,\eta \in H$.

On appelle \emph{alg{\`e}bre de von Neumann} toute sous-alg{\`e}bre involutive $M \subset \B(H)$ qui est {\'e}gale {\`a} son bicommutant~: $M = M\dpr$, ce qui
{\'e}quivaut {\`a} dire que $M$ est faiblement ferm{\'e} et $1 \in M$. Un
\emph{facteur} est une alg{\`e}bre de von Neumann dont le centre est r{\'e}duit aux
scalaires. Si $G$ est un groupe localement compact, l'alg{\`e}bre de von
Neumann du groupe $G$ not{\'e}e $L(G)$ est le bicommutant $\{\la_g \mid g \in G\}\dpr$
o{\`u} $(\la_g)$ est la repr{\'e}sentation r{\'e}guli{\`e}re du groupe $G$ sur
l'espace de Hilbert $L^2(G)$.

Murray et von Neumann ont classifi{\'e} les facteurs en type I, II et III. Les \emph{facteurs de type I} sont ceux qui poss{\`e}dent des projecteurs
minimaux. Ils sont isomorphes {\`a} $M_n(\C)$ (type I$_n$) ou $\B(\ell^2)$ (type I$_\infty$). Les \emph{facteurs de type II$_1$} sont ceux qui admettent
une trace finie et qui sont de dimension infinie (pour exclure le cas I$_n$). L'exemple type d'un facteur II$_1$ est donn{\'e} par l'alg{\`e}bre de von
Neumann $L(G)$ d'un groupe discret $G$ dont $\{e\}$ est la seule classe de conjugaison finie (on dit que $G$ est CCI). Les groupes libres $\F_n$ {\`a}
$n$ g{\'e}n{\'e}rateurs sont des exemples de groupes CCI ($n$ peut {\^e}tre $\infty$). Les \emph{facteurs de type II$_\infty$} sont ceux qui sont de la forme $N
\ot \B(\ell^2)$ avec $N$ un facteur II$_1$. Ce sont exactement les facteurs qui admettent une trace infinie semi-finie et qui ne sont pas de type I.
Un facteur de type I ou II est dit \emph{semi-fini}. Il admet toujours une trace semi-finie. Finalement les \emph{facteurs de type III} sont ceux qui
n'admettent pas de trace non-nulle.

\begin{rema}
Dans tout l'expos{\'e} les espaces de Hilbert sont suppos{\'e}s
\emph{s{\'e}parables} et les alg{\`e}bres de von Neumann
\emph{{\`a} pr{\'e}dual s{\'e}parable}, i.e.\ admettant une
repr{\'e}sentation fid{\`e}le sur un espace de Hilbert s{\'e}parable.
\end{rema}

\subsection{Classification des facteurs de type III d'apr{\`e}s Connes}

Un \emph{{\'e}tat normal} d'une alg{\`e}bre de von Neumann $M$ est une forme lin{\'e}aire faiblement continue $\om : M \recht \C$ positive ($\om(x) \geq 0$
quand $x \geq 0$) qui satisfait $\om(1) = 1$. Un {\'e}tat est dit \emph{fid{\`e}le} si $x=0$ d{\`e}s que $\om(x) = 0$ et $x \geq 0$. Toute alg{\`e}bre de von
Neumann {\`a} pr{\'e}dual s{\'e}parable admet un {\'e}tat fid{\`e}le.

La th{\'e}orie de Tomita-Takesaki associe {\`a} tout {\'e}tat fid{\`e}le $\om$ un groupe {\`a} un
param{\`e}tre $(\si^\om_t)$ d'automorphismes de $M$ appel{\'e} \emph{groupe
  modulaire} et caract{\'e}ris{\'e} par
\begin{itemize}
\item $\om \si^\om_t = \om$ pour tout $t \in \R$,
\item $\om$ satisfait la condition KMS par rapport {\`a} $(\si^\om_t)$~: pour tout $x,y \in M$, il existe une fonction continue $f : \{z \in \C \mid 0 \leq
\operatorname{Im} z \leq 1 \} \recht \C$ qui est analytique {\`a}
l'int{\'e}rieur de la bande et qui satisfait
$$f(t) = \om(x \si^\om_t(y)) \quad\text{et}\quad f(t+i) = \om(\si^\om_t(y) x) \quad\text{pour tout}\; t \in \R \; .$$
\end{itemize}
Une bonne introduction {\`a} la \emph{th{\'e}orie
  modulaire de Tomita-Takesaki} se trouve dans \cite{stratila-zsido}.

Le \emph{th{\'e}or{\`e}me de Radon-Nikodym} de Connes \cite{Connes-class}
permet de comparer les groupes modulaires de deux {\'e}tats fid{\`e}les $\om$,
$\mu$ sur $M$. En effet, il existe une
application faiblement continue $t \mapsto u_t$ de $\R$ dans le groupe unitaire de $M$ telle que
\begin{itemize}
\item $u_{t+s} = u_t \si^\om_t(u_s)$,
\item $\si^\mu_t(x) = u_t \si^\om_t(x) u_t^*$.
\end{itemize}
Les groupes modulaires de deux {\'e}tats fid{\`e}les diff{\`e}rent donc
par une perturbation int{\'e}rieure. Il s'en suit qu'une alg{\`e}bre de von
Neumann a une dynamique intrins{\`e}que donn{\'e}e par les groupes modulaires
des {\'e}tats fid{\`e}les et d{\'e}termin{\'e}e {\`a} perturbation int{\'e}rieure pr{\`e}s.

D{\'e}finissons le groupe polonais $\Aut M$ des automorphismes de $M$ muni de la topologie induite par les distances $d(\al,\be) = \|\om \al - \om \be\|$
et $d(\al,\be) = \|\om \al^{-1} - \om \be^{-1}\|$, o{\`u} $\om$ parcourt les {\'e}tats de $M$. A chaque unitaire $u \in M$, on associe \emph{l'automorphisme
int{\'e}rieur} $\Ad u$ d{\'e}fini par $(\Ad u)(x) = uxu^*$. Les automorphismes int{\'e}rieurs forment un sous-groupe distingu{\'e} $\Int M$ de $\Aut M$. Le groupe
quotient est not{\'e} $\Out M$. Le th{\'e}or{\`e}me de Radon-Nikodym permet de d{\'e}finir un homomorphisme $\delta : \R \recht \Out M$ qui envoie $t$ {\`a} la classe de
$\si^\om_t$ et qui ne d{\'e}pend pas du choix de $\om$.

\subsubsection*{Invariant $T$}

Soit $\om$ un {\'e}tat fid{\`e}le sur un facteur $M$ avec groupe modulaire $(\si_t)$. Dans \cite{Connes-class} Connes a introduit le sous-groupe $T(M)$ de
$\R$~:
$$T(M) = \{ t \in \R \mid \si^\om_t \in \Int M \} \; .$$
D'apr{\`e}s le th{\'e}or{\`e}me de Radon-Nikodym, l'invariant $T(M)$ ne d{\'e}pend pas du choix de l'{\'e}tat $\om$. C'est
donc un invariant de l'alg{\`e}bre de von Neumann $M$.

\subsubsection*{Flot des poids}

A chaque facteur $M$ est associ{\'e} une alg{\`e}bre de von
Neumann semi-finie~: c'est le produit crois{\'e} $M \rtimes_{(\si_t)}
\R$ de $M$ par le groupe modulaire $(\si_t)$ d'un {\'e}tat fid{\`e}le sur
$M$. Sur le produit crois{\'e} $M \rtimes_{(\si_t)}
\R$, il existe une trace semi-finie canonique. Le produit crois{\'e} $M \rtimes_{(\si_t)}
\R$ admet une \emph{action duale} $(\te_s)$ de $\R^*_+$ par
automorphismes. La restriction de l'action $(\te_s)$ au centre du
produit crois{\'e} s'appelle le \emph{flot des poids} de $M$.
Gr{\^a}ce au th{\'e}or{\`e}me de Radon-Nikodym, le flot des poids ne d{\'e}pend pas
du choix de l'{\'e}tat fid{\`e}le.

\subsubsection*{Facteurs de type III$_\la$}

Notons que le flot des poids est une action
ergodique de $\R^*_+$ sur un espace mesur{\'e}. Or une telle action est ou
bien transitive ou bien proprement ergodique, d'o{\`u} la classification suivante~:
\begin{itemize}
\item $M$ est semi-fini si c'est l'action de $\R^*_+$ sur $\R^*_+$,
\item $M$ est de type III$_\la$ avec $0 < \la < 1$ si c'est l'action
  de $\R^*_+$ sur $\R^*_+/\la^\Z$,
\item $M$ est de type III$_1$ si c'est l'action de $\R^*_+$ sur un point,
\item $M$ est de type III$_0$ si c'est une action proprement
  ergodique.
\end{itemize}
Combes a donn{\'e} un aper{\c c}u des r{\'e}sultats de classification de Connes dans \cite{combes}.

Murray et von Neumann ont construit deux facteurs de type II$_1$
non-isomorphes~: le facteur hyperfini $\cR$ et le facteur $L(\F_2)$ du
groupe libre {\`a} 2 g{\'e}n{\'e}rateurs. Ils sont non-isomorphes car
$\Int L(\F_2)$ est ferm{\'e} dans $\Aut L(\F_2)$ tandis que $\Int \cR$ est
un sous-groupe dense non-trivial de $\Aut \cR$.

\subsubsection*{Facteurs pleins et invariant $\tau$}

Un facteur $M$ est dit \emph{plein} si $\Int M$
est un sous-groupe ferm{\'e} de $\Aut M$. Pour un facteur plein, le groupe
quotient $\Out M$ est un groupe polonais.
Pour un tel facteur plein Connes \cite{Connes-periodic} introduit un nouvel
invariant~:
$$\tau(M) = \text{la topologie la plus faible sur $\R$ qui rend continue
  l'application $\delta : \R \recht \Out M$}$$
Les facteurs d'Araki-Woods libres sont des facteurs pleins. Nous
remarquons qu'un facteur plein ne peut {\^e}tre de type III$_0$.

\subsection{Probabilit{\'e}s libres d'apr{\`e}s Voiculescu}

Une introduction plus compl{\`e}te aux probabilit{\'e}s libres de Voiculescu se trouve dans le livre \cite{Voic-livre} ou dans \cite{skandal-bourb}.

Un \emph{espace de probabilit{\'e}s non-commutatif}  est une paire
$(A,\vfi)$ o{\`u} $A$ est une alg{\`e}bre unif{\`e}re et $\vfi$ est une forme
lin{\'e}aire v{\'e}rifiant $\vfi(1)= 1$. Dans cet
expos{\'e} on s'int{\'e}resse surtout aux alg{\`e}bres de von Neumann~:
$A$ est une alg{\`e}bre de von Neumann et $\vfi$ un {\'e}tat normal. Les {\'e}l{\'e}ments
de $A$ s'appellent toujours variables al{\'e}atoires. La
\emph{distribution} d'un {\'e}l{\'e}ment $x \in A$ est l'application qui
{\`a} un polyn{\^o}me $P \in \C[X]$ associe $\vfi(P(x))$. Si $A$ est une
alg{\`e}bre de von Neumann et $x$ un {\'e}l{\'e}ment auto-adjoint, la distribution
de $x$ est une mesure de probabilit{\'e}s dont le support est contenu dans
$[-\|x\|,\|x\|]$.

\begin{nota}
Si $(M,\vfi)$ et $(N,\mu)$ sont des alg{\`e}bres de von Neumann munies
d'{\'e}tats $\vfi$ et $\mu$, la notation $(M,\vfi) \cong (N,\mu)$ signifie
qu'il existe un $^*$-isomorphisme $\al : M \recht N$ tel que $\mu \al
= \vfi$.
\end{nota}

\begin{defi}
Soit $(A,\vfi)$ un espace de probabilit{\'e}s non-commutatif. Une famille
$(A_i)_{i \in I}$ de sous-alg{\`e}bres est dite \emph{libre} si pour tout
$k$ et toute suite d'{\'e}l{\'e}ments $a_j \in A_{i_j}$ ($j=1,\ldots,k$)
satisfaisant $\vfi(a_j)=0$ et $i_j \neq i_{j+1}$, on a $\vfi(a_1\cdots
a_k) = 0$.

Une famille d'{\'e}l{\'e}ments $(x_i)_{i \in I}$ est dite \emph{libre} (resp.\
\emph{$^*$-libre}) si les alg{\`e}bres (resp.\ $^*$-alg{\`e}bres) $A_i$ engendr{\'e}es
par $x_i$ forment une famille libre de de sous-alg{\`e}bres de $A$.
\end{defi}

Les produits libres fournissent des exemples de familles libres.

\begin{prop}
Soit $(M_i,\vfi_i)$ une famille d'alg{\`e}bres de von Neumann munies d'un
{\'e}tat fid{\`e}le. Alors, il existe, {\`a} isomorphisme pr{\`e}s,
une unique paire $(M,\vfi)$ d'une alg{\`e}bre von Neumann munie d'un {\'e}tat
fid{\`e}le, telle que
\begin{itemize}
\item $(M_i,\vfi_i)$ se plonge dans $(M,\vfi)$ en pr{\'e}servant l'{\'e}tat,
\item $M$ est engendr{\'e}e par la famille de sous-alg{\`e}bres $(M_i)$ qui
  est une famille libre dans $(M,\vfi)$.
\end{itemize}
On appelle $(M,\vfi)$ le \emph{produit libre des $(M_i,\vfi_i)$} et on note $(M,\vfi) = \freep{i \in I} (M_i,\vfi_i)$.
\end{prop}

\subsection{Repr{\'e}sentations orthogonales de $\R$} \label{sec.rep}

La donn{\'e}e de la construction des {\'e}tats quasi-libres libres et
des facteurs d'Araki-Woods libres associ{\'e}s, est une repr{\'e}sentation de
$\R$ par transformations orthogonales.

\begin{term}
On appelle \emph{repr{\'e}sentation orthogonale de $\R$} tout
groupe {\`a} un param{\`e}tre de transformations orthogonales d'un espace de Hilbert r{\'e}el $H_\R$.
\end{term}

Soit $(U_t)$ une repr{\'e}sentation orthogonale de $\R$ sur l'espace de Hilbert r{\'e}el $H_\R$. Le complexifi{\'e} $H = H_\R \ot \C$ admet une
involution anti-unitaire $J$ (l'op{\'e}rateur de conjugaison complexe) et les transformations orthogonales $(U_t)$ s'{\'e}tendent en un groupe {\`a} un param{\`e}tre d'unitaires sur $H$, qu'on notera
toujours $(U_t)$.

Il existe alors un unique op{\'e}rateur auto-adjoint strictement positif
$A$ sur $H$ tel que $U_t = A^{it}$ pour tout $t \in \R$. On a $JAJ =
A^{-1}$. L'op{\'e}rateur $A$ permet de d{\'e}finir un nouveau plongement
isom{\'e}trique
$$H_\R \hookrightarrow H : \xi \mapsto \bigl( \frac{2}{A^{-1} + 1}
\bigr)^{1/2} \xi \; .$$
En effet, si $\xi \in H_\R$, on a $J \xi = \xi$ et donc
\begin{align*}
\bigl\| \bigl(\frac{2}{A^{-1} + 1} \bigr)^{1/2} \xi \bigr\|^2 &= \langle \frac{1}{A^{-1} + 1} \xi,\xi \rangle + \langle \frac{1}{A^{-1} + 1}
J\xi,J\xi \rangle = \langle \frac{A}{A+1} \xi,\xi \rangle + \langle J \frac{1}{A + 1} \xi,J\xi \rangle
 \\ &= \langle \frac{A+1}{A+1} \xi,\xi \rangle = \|\xi\|^2 \; .
\end{align*}
On notera $K_\R$ l'image de $H_\R$ par ce plongement. Alors,
$K_\R$ est un espace de Hilbert r{\'e}el, isom{\'e}triquement plong{\'e} dans un
espace de Hilbert complexe
$H$ v{\'e}rifiant la propri{\'e}t{\'e} suivante~:
\begin{enumerate}
\item[($\star$)] $K_\R \cap i K_\R = \{0\}$ et $K_\R + i K_\R$ est dense dans $H$.
\end{enumerate}
Dans \cite{rieff-vandaele} on d{\'e}montre que chaque plongement isom{\'e}trique $K_\R \subset H$ satisfaisant la condition ($\star$) provient d'une
repr{\'e}sentation orthogonale de $\R$ sur un espace de Hilbert r{\'e}el $H_\R$ par la construction pr{\'e}sent{\'e}e ci-dessus.

{\'E}crivons $T= J A^{-1/2}$. Alors $T$ est un op{\'e}rateur anti-lin{\'e}aire ferm{\'e} et inversible sur $H$ qui satisfait $T = T^{-1}$. Un tel op{\'e}rateur s'appelle
une \emph{involution sur $H$}. R{\'e}ciproquement une telle involution $T$ admet une d{\'e}composition polaire $T = JA^{-1/2}$ dans laquelle $J$ est une
involution anti-unitaire sur $H$ et $A$ est un op{\'e}rateur auto-adjoint strictement positif satisfaisant $JAJ = A^{-1}$. Posons $H_\R = \{\xi \in H
\mid J\xi = \xi \}$ et $U_t = A^{it}$. On obtient ainsi une repr{\'e}sentation orthogonale de $\R$. On remarquera que l'espace $K_\R$ correspondant
consiste des vecteurs $\xi$ dans le domaine de $T$ qui satisfont $T \xi = \xi$.

On a alors obtenu plusieurs points de vue diff{\'e}rents sur les
repr{\'e}sentations orthogonales
de $\R$.
\begin{enumerate}
\item Un groupe {\`a} un param{\`e}tre de transformations orthogonales d'un
  espace de Hilbert r{\'e}el.
\item Un plongement isom{\'e}trique d'un espace de Hilbert r{\'e}el dans un espace de
  Hilbert complexe v{\'e}rifiant ($\star$).
\item Une involution $T$ sur un espace de Hilbert.
\end{enumerate}
Finalement on peut consid{\'e}rer la d{\'e}composition spectrale de
l'op{\'e}rateur $\log A$. Comme $J (\log A) J = - \log A$, la classe de la
mesure spectrale de $\log A$ est sym{\'e}trique. Les repr{\'e}\-sen\-tations orthogonales
de $\R$ sont donc classifi{\'e}es par une classe de mesures sym{\'e}trique sur
$\R$ et une fonction de multiplicit{\'e} sym{\'e}trique.

\section{Facteurs d'Araki-Woods libres} \label{sec.arakiwoodslibre}

Le foncteur $\car$ associe {\`a} tout espace de Hilbert $H$ la
\cst-alg{\`e}bre $\car(H)$ (voir introduction). Oubliant la structure
complexe de $H$ on peut {\'e}crire $\car(H)$ comme une alg{\`e}bre de
Clifford. Soit $H_\R$ un espace de Hilbert r{\'e}el. On note
$\cliff(H_\R)$ et on appelle \emph{alg{\`e}bre de Clifford} la \cst-alg{\`e}bre universelle engendr{\'e}e par la famille
$\{s(\xi) \mid \xi \in H_\R\}$ telle que $s(\xi)$ est auto-adjoint
pour tout $\xi \in H_\R$, $\xi \mapsto s(\xi)$ est $\R$-lin{\'e}aire et
$$s(\xi) s(\eta) + s(\eta) s(\xi) = 2 \langle \xi,\eta \rangle 1 \;
.$$
Cette derni{\`e}re condition {\'e}tant {\'e}quivalente {\`a} $s(\xi)^2 = \|\xi\|^2 1$
pour tout $\xi \in H_\R$, on voit comment le foncteur $H_\R \mapsto
\libre(H_\R)$ est une version libre du foncteur $\cliff$.

A chaque plongement isom{\'e}trique $H_\R \hookrightarrow H$ de $H_\R$
dans un espace de Hilbert complexe $H$ est associ{\'e}e une repr{\'e}sentation
de $\cliff(H_\R)$ sur \emph{l'espace de Fock anti-sym{\'e}trique (ou fermionique)}~:
$$\cF_{\text{\rm as}}(H) = \C \Om \oplus \bigoplus_{n=1}^\infty
H^{\wedge n}$$
posant $s(\xi) = a(\xi)^* + a(\xi)$ o{\`u} $a(\xi)$ est l'op{\'e}rateur de
cr{\'e}ation {\`a} gauche. Remarquons que cette repr{\'e}sentation de
$\cliff(H_\R)$ est en fait la repr{\'e}sentation GNS d'un {\'e}tat
quasi-libre. L'alg{\`e}bre de von Neumann engendr{\'e}e par les
op{\'e}rateurs $s(\xi), \xi \in H_\R$ est un \emph{facteur d'Araki-Woods}.

Shlyakhtenko donne une version libre de la construction pr{\'e}c{\'e}dente et
appelle le facteur engendr{\'e} \emph{facteur d'Araki-Woods libre}.

\subsection{{\'E}tats quasi-libres libres}

Donnons nous une repr{\'e}sentation orthogonale $(U_t)$ de $\R$
sur l'espace de Hilbert r{\'e}el $H_\R$. Comme au
\S\ref{sec.rep}, nous regardons le complexifi{\'e} $H$ de $H_\R$ avec l'involution anti-unitaire $J$ et l'op{\'e}rateur auto-adjoint strictement positif $A$
tel que $U_t = A^{it}$. Introduisons \emph{l'espace de Fock
  plein} de $H$~:
$$\cF(H) = \C \Om \oplus \bigoplus_{n=1}^\infty H^{\ot n} \; .$$
Le vecteur unit{\'e} $\Om$ s'appelle \emph{vecteur du vide}. Pour chaque
vecteur $\xi \in H$, nous disposons de \emph{l'op{\'e}rateur de cr{\'e}ation {\`a}
gauche}
$$\ell(\xi) : \cF(H) \recht \cF(H) : \begin{cases} \ell(\xi)\Om = \xi
  \; , \\ \ell(\xi)(\xi_1 \ot \cdots \ot \xi_n) = \xi \ot \xi_1 \ot
  \cdots \ot \xi_n \; . \end{cases}$$
L'adjoint $\ell(\xi)^*$ s'appelle \emph{op{\'e}rateur d'annihilation}.

Pour chaque vecteur $\xi \in H$, notons $s(\xi)$ la partie r{\'e}elle
de $\ell(\xi)$ donn{\'e}e par
$$s(\xi) = \frac{\ell(\xi) + \ell(\xi)^*}{2} \; .$$
Un r{\'e}sultat crucial de Voiculescu \cite{Voic-livre} dit que la distribution
de l'op{\'e}rateur $s(\xi)$ par rapport {\`a} l'{\'e}tat vectoriel du vide donn{\'e} par
$\vfi(x) = \langle x \Om,\Om \rangle$ est la loi semi-circulaire de
Wigner support{\'e}e par l'intervalle $[-\|\xi\|,\|\xi\|]$.

Rappelons que l'op{\'e}rateur $A$ permet de d{\'e}finir un plongement de $H_\R$ dans $H$ dont l'image est not{\'e}e $K_\R$. On peut alors formuler la
d{\'e}finition centrale de cet expos{\'e} \cite{shl1}.

\begin{defi}
Soit $(U_t)$ une repr{\'e}sentation orthogonale de $\R$ sur l'espace de
Hilbert r{\'e}el $H_\R$.
Le \emph{facteur d'Araki-Woods libre}\footnote{Nous verrons que
  $\free{H_\R,U_t}$ est effectivement un facteur d{\`e}s que $\dim H_\R
  \geq 2$.} not{\'e} $\free{H_\R,U_t}$ est
d{\'e}fini par
$$\free{H_\R,U_t} = \{s(\xi) \mid \xi \in K_\R \}\dpr \; .$$
L'{\'e}tat vectoriel
$\vfi_U(x) = \langle x \Om,\Om \rangle$
est appel{\'e} \emph{{\'e}tat quasi-libre libre}.
\end{defi}

Rappelons que $T= JA^{-1/2}$ est l'involution sur $H$ associ{\'e}e {\`a} $(U_t)$. Pour $\xi,\eta \in K_\R$, on v{\'e}rifie que
$$2 s(\xi) + 2i s(\eta) = \ell(\zeta) + \ell(T\zeta)^*$$
o{\`u} $\zeta= \xi + i \eta$. On conclut que $\free{H_\R,U_t}$ est
{\'e}galement l'alg{\`e}bre de von Neumann engendr{\'e}e par les op{\'e}rateurs
$\ell(\zeta) + \ell(T\zeta)^*$ o{\`u} $\zeta$ appartient au domaine de
$T$.

Le r{\'e}sultat suivant est facile {\`a} d{\'e}montrer.
\begin{prop} \label{prop.groupemodulaire}
L'{\'e}tat quasi-libre libre $\vfi_U$ sur $\free{H_\R,U_t}$ est fid{\`e}le. Le
groupe modulaire $(\si_t)$ de l'{\'e}tat $\vfi_U$ est donn{\'e} par
$$\si_t(s(\xi)) = s(U_t \xi) \quad\text{pour tout}\;\; t \in
\R, \xi \in K_\R \; .$$
\end{prop}

La construction des facteurs d'Araki-Woods libres est fonctorielle dans un sens pr{\'e}cis. En effet on consid{\`e}re la cat{\'e}gorie dont les objets sont les
paires $(H_\R,U_t)$ et les morphismes sont les contractions entre espaces de Hilbert qui entrelacent les repr{\'e}sentations. A chaque morphisme
$(H^{(1)}_\R,U^{(1)}_t) \recht (H^{(2)}_\R,U^{(2)}_t)$ correspond une application compl{\`e}tement positive $\free{H^{(1)}_\R,U^{(1)}_t} \recht
\free{H^{(2)}_\R,U^{(2)}_t}$ normale et unif{\`e}re, pr{\'e}servant les {\'e}tats quasi-libres libres. On notera $\Gamma\dpr$ ce foncteur.

La cat{\'e}gorie des paires $(H_\R,U_t)$ admet une structure additive~: la somme directe. Shlyakhtenko d{\'e}montre que le
foncteur $\Gamma\dpr$ entrelace les op{\'e}rations somme directe et produit
libre.
\begin{prop} \label{prop.sommedirecte}
Soit $(H^{(i)}_\R,U^{(i)}_t)_{i \in I}$ une repr{\'e}sentation orthogonale de $\R$. Posons $(H_\R,U_t) = \bigoplus_{i} (H^{(i)}_\R,U^{(i)}_t)$. Alors,
$$(\free{H_\R,U_t},\vfi_U) \cong \freep{i \in I}
\bigl(\free{H^{(i)}_\R,U^{(i)}_t},\vfi_{U^{(i)}} \bigr) \; .$$
\end{prop}

\subsection{Variables circulaires g{\'e}n{\'e}ralis{\'e}es} \label{sec.circulaire}

Pour comprendre la structure des alg{\`e}bres de von Neumann $(\free{H_\R,U_t},\vfi_U)$ il
est naturel de consid{\'e}rer d'abord les repr{\'e}sentations orthogonales
irr{\'e}ductibles de $\R$. Le cas $H_\R = \R$ et $U_t = \id$ est facile~:
l'alg{\`e}bre est engendr{\'e} par un seul op{\'e}rateur dont la distribution par
rapport {\`a} $\vfi_U$ est la loi semi-circulaire, d'apr{\`e}s le r{\'e}sultat de Voiculescu. On trouve donc
$$(\free{\R,\id},\vfi_U) \cong (L^\infty[-1,1],\mu)$$
o{\`u} $\mu$ est la mesure semi-circulaire sur $[-1,1]$. Si on combine ce
r{\'e}sultat avec la proposition \ref{prop.sommedirecte}, on obtient
\begin{equation}\label{eq.trivial}
(\free{H_\R,\id},\vfi_U) \cong (L(\F_n),\trace)
\end{equation}
o{\`u} $L(\F_n)$ est l'alg{\`e}bre de von Neumann du groupe libre {\`a} $n = \dim H_\R$
g{\'e}n{\'e}rateurs.

Prenons maintenant $H_\R = \R^2$ et $0 < \la < 1$. Posons
\begin{equation}\label{eq.matrice}
U_t =
\begin{pmatrix} \cos(t \log \la) & - \sin (t \log \la) \\
  \sin(t \log \la) & \cos(t \log \la) \end{pmatrix} \; .
\end{equation}
En prenant la base orthonormale $\xi_1=\frac{1}{\sqrt{2}}(1,-i)$,
$\xi_2=\frac{1}{\sqrt{2}}(1,i)$ du complexifi{\'e} $H = \C^2$, on voit que
l'alg{\`e}bre de von
Neumann $\free{H_\R,U_t}$ est engendr{\'e}e par l'op{\'e}rateur $\ell(\xi_2) + \sqrt{\la} \ell(\xi_1)^*$ sur l'espace de Fock plein $\cF(\C^2)$.

\begin{nota}
On notera $(T_\la,\vfi_\la) := \free{H_\R,U_t}$ o{\`u} $H_\R= \R^2$ et
$U_t$ est donn{\'e} par l'{\'e}galit{\'e} \eqref{eq.matrice}.
\end{nota}

Pour comprendre l'alg{\`e}bre $T_\la$, il faut {\'e}tudier la
$^*$-distribution de l'{\'e}l{\'e}ment $\ell(\xi_2) + \sqrt{\la}
\ell(\xi_1)^*$ par rapport {\`a} l'{\'e}tat vectoriel du vide. Un tel {\'e}l{\'e}ment
s'appelle \emph{{\'e}l{\'e}ment circulaire g{\'e}n{\'e}ralis{\'e}}. Dans le cas $\la = 1$,
on retrouve l'{\'e}l{\'e}ment circulaire $y$ de Voiculescu \cite{Voic-livre}. Voiculescu
a d{\'e}montr{\'e} que la d{\'e}composition polaire $y = u b$ d'un {\'e}l{\'e}ment
circulaire donne un \emph{unitaire de Haar} $u$ et un op{\'e}rateur $b$
\emph{quart-circulaire}. Ceci veut dire que la distribution de $u$ est
la distribution uniforme sur le cercle et que la distribution de $b$
suit la loi quart-circulaire support{\'e}e par l'intervalle
$[0,1]$. Shlyakthenko a d{\'e}montr{\'e} dans \cite{shl1} un r{\'e}sultat analogue pour les {\'e}l{\'e}ments
circulaires g{\'e}n{\'e}ralis{\'e}s. Ce r{\'e}sultat permet de donner une description
alternative de $(T_\la,\vfi_\la)$.

\begin{theo} \label{theo.Tla}
Soit $0 < \la < 1$ et soit $y = \ell(\xi_1) + \sqrt{\la}\ell(\xi_2)^*$
l'{\'e}l{\'e}ment circulaire g{\'e}n{\'e}ralis{\'e} associ{\'e} dans $(T_\la,\vfi_\la)$. Notons $y = vb$ la
d{\'e}composition polaire de $y$. Alors $v$ est une isom{\'e}trie non-unitaire
qui satisfait $\vfi_\la(v^k (v^*)^l) = \delta_{kl} \la^k$. La
distribution de l'op{\'e}rateur $b$ est sans atomes.
Les {\'e}l{\'e}ments $u$ et $b$ sont $^*$-libres.
\end{theo}

Un corollaire imm{\'e}diat de ce r{\'e}sultat est que
\begin{equation}\label{eq.freeproduct}
(T_\la,\vfi_\la) \cong
(\B(\ell^2(\N)),\om_\la) \ast (L^\infty[-1,1],\mu)
\end{equation}
o{\`u}
$\om_\la(e_{ij}) = \delta_{ij} \la^j (1-\la)$ et $\mu$ est la loi
semi-circulaire sur $[-1,1]$. Bien {\'e}videmment, au lieu de $\mu$ on pourrait prendre
n'importe quelle autre mesure de probabilit{\'e}s sans atomes.

L'isomorphisme \eqref{eq.freeproduct} est crucial. Il permet de r{\'e}aliser
$(T_\la,\vfi_\la)$ en repr{\'e}sentant d'une mani{\`e}re libre $(\B(\ell^2(\N)),\om_\la)$ et
$(L^\infty[-1,1],\mu)$ dans un espace de probabilit{\'e}
non-commu\-ta\-tif. Shlyakhtenko trouve dans \cite{shl1} de telles repr{\'e}sentations qui
permettent de comprendre la r{\'e}duc\-tion de l'alg{\`e}bre $T_\la$ par un
projecteur minimal de $\B(\ell^2(\N))$. On les appelle \emph{mod{\`e}les
  matriciels}. C'est un outil puissant qui permet de d{\'e}montrer des
r{\'e}sultats \emph{d'absorption libre}.

\begin{theo} \label{theo.absorption}
On a
$$
(T_\la,\vfi_\la) \cong (T_\la,\vfi_\la) \ast (L^\infty[-1,1],\mu)
\cong (T_\la,\vfi_\la) \ast (L(\F_\infty),\trace) \; ,
$$
o{\`u} $\mu$ est la mesure semi-circulaire et $\trace$ est la trace sur le
facteur $L(\F_\infty)$ du groupe libre {\`a} une infinit{\'e} de g{\'e}n{\'e}rateurs.
\end{theo}

\subsection{Type des facteurs d'Araki-Woods libres} \label{subsec.types}

A l'aide du th{\'e}or{\`e}me \ref{theo.absorption}, on peut finalement d{\'e}montrer que $\free{H_\R,U_t}$ est toujours un facteur quand la dimension de $H_\R$
est au moins $2$. On peut en m{\^e}me temps d{\'e}terminer le type de ce facteur et son invariant $\tau$.

\begin{theo} \label{theo.types}
Soit $(U_t)$ une repr{\'e}sentation orthogonale de $\R$ sur l'espace de Hilbert r{\'e}el $H_\R$ de dimension au moins $2$. Notons
$M = \free{H_\R,U_t}$.
\begin{enumerate}
\item $M$ est un facteur plein.
\item $M$ est de type II$_1$ ssi $U_t = \id$ pour tout $t \in \R$.
\item $M$ est de type III$_\la$ $(0 < \la < 1)$ ssi $(U_t)$ est
  p{\'e}riodique de p{\'e}riode $\frac{2\pi}{|\log \la|}$.
\item $M$ est de type III$_1$ dans les autres cas.
\item L'invariant $\tau(M)$ est la topologie la plus faible sur $\R$
  qui rend continue l'application $t \mapsto U_t$ de $\R$ dans le
  groupe orthogonal de $H_\R$ muni de la topologie faible.
\item Le facteur $M$ admet des {\'e}tats presque-p{\'e}riodiques ssi $(U_t)$ est presque-p{\'e}riodique.
\end{enumerate}
\end{theo}

Nous donnons ici plus de d{\'e}tails pour la d{\'e}monstration de ce th{\'e}or{\`e}me. Shlyakh\-tenko  d{\'e}termine l'invariant $\tau(M)$ dans \cite{shl3}, mais en
supposant que la repr{\'e}sentation orthogonale $(U_t)$ contient ou bien une repr{\'e}sentation p{\'e}riodique ou bien une repr{\'e}sentation triviale de dimension
$2$. Nous suivons la m{\^e}me m{\'e}thode que Shlyakhtenko, mais utilisons le nouveau lemme \ref{lemm.technique} qui est plus fort que le lemme des 14$\eps$
de Barnett \cite{Bar} utilis{\'e} par Shlyakhtenko. Shlyakhtenko d{\'e}montre
(1)--(4) dans \cite{shl1} pour les repr{\'e}sentations
presque-p{\'e}riodiques et dans \cite{shl2} pour le cas g{\'e}n{\'e}ral, mais par d'autres m{\'e}thodes que nous.

\begin{proof}[Preuve du th{\'e}or{\`e}me \ref{theo.types}]
Il suffit de d{\'e}montrer (1) et (5). En effet, un facteur plein est semi-fini (c'est-{\`a}-dire, de type I ou II) ssi $\tau(M)$ est la topologie grossi{\`e}re.
Dans ce cas-l{\`a}, on conclut de (5) que $U_t = \id$ pour tout $t \in \R$ et d'apr{\`e}s l'isomorphisme \eqref{eq.trivial}, $M$ est un facteur II$_1$. Ceci
d{\'e}montre (2). Un facteur plein n'est jamais de type III$_0$. Comme $\tau(M)$ est la topologie la plus faible qui rend continue l'application $\delta
: \R \recht \operatorname{Out}(M)$, on conclut de (5) que $\delta(t)=1$ ssi $U_t = \id$. Ceci d{\'e}montre (3) et (4). Finalement, d{\'e}montrons (6). Si $M$
admet un {\'e}tat presque-p{\'e}riodique, le groupe $\R$ munie de la topologie $\tau(M)$ peut {\^e}tre compl{\'e}t{\'e} en un groupe compact. Il existe donc un groupe
compact $G$, un plongement $\R \subset G$ et une extension de $t \mapsto U_t$ en un homomorphisme continu $G \recht O(H_\R)$. Ceci veut dire que
$(U_t)$ est presque-p{\'e}riodique \cite{dixmier}. R{\'e}ciproquement, si $(U_t)$ est presque-p{\'e}riodique, l'{\'e}tat quasi-libre libre est un {\'e}tat
presque-p{\'e}riodique.

Il nous reste {\`a} d{\'e}montrer (1) et (5). Ceci est {\'e}vident quand $U_t = \id$ pour tout $t \in \R$. Le deuxi{\`e}me cas qu'on consid{\`e}re est celui o{\`u} $(U_t)$
contient la repr{\'e}sentation donn{\'e}e par l'{\'e}galit{\'e} \eqref{eq.matrice} avec $0 < \la < 1$. Notons son compl{\'e}ment par $(U'_t)$ agissant sur $H'_\R$.
D'apr{\`e}s le th{\'e}or{\`e}me \ref{theo.absorption} on a
\begin{align*}
(M,\vfi) &\cong (T_\la,\vfi_\la) \ast (\free{H'_\R,U'_t},\vfi_{U'})
\\ &\cong
\bigl( (T_\la,\vfi_\la) \ast (L^\infty[-1,1],\mu) \bigr) \ast \bigl(
L^\infty([-1,1],\mu) \ast (\free{H'_\R,U'_t},\vfi_{U'})\bigr) \; .
\end{align*}
Comme $(L^\infty[-1,1],\mu)$ contient un unitaire de Haar, on peut appliquer la proposition \ref{prop.invarianttau}. On conclut que $M$ est un
facteur plein et que l'invariant $\tau(M)$ est la topologie la plus faible sur $\R$ qui rend continues les deux applications $t \mapsto
\si_t^{\vfi_\la}$ et $t \mapsto \si_t^{\vfi_{U'}}$. Par la proposition \ref{prop.groupemodulaire} ceci est exactement la topologie la plus faible sur
$\R$ qui rend continue l'application $t \mapsto U_t$.

Finalement, nous consid{\'e}rons le cas o{\`u} la repr{\'e}sentation $(U_t)$ ne contient pas de repr{\'e}\-sentation p{\'e}riodique et n'est pas triviale. Il est alors
clair qu'on peut d{\'e}composer $U_t$ en trois composantes non-triviales $U_t = U^{(1)}_t \oplus U^{(2)}_t \oplus U^{(3)}_t$. Les {\'e}nonc{\'e}s (1) et (5)
d{\'e}coulent des propositions \ref{prop.sommedirecte} et
\ref{prop.invarianttau} ainsi que du lemme \ref{lemm.approx}.
\end{proof}

Pour chaque repr{\'e}sentation orthogonale non-p{\'e}riodique $(U_t)$ de $\R$, le facteur
d'Araki-Woods libre $\free{H_\R,U_t}$ est donc un facteur
de type III$_1$ dont l'invariant $\tau$ est la topologie la plus
faible qui rend continue l'application $t \mapsto U_t$.

\begin{rema}
Dans
\cite{Connes-periodic} Connes part d'une mesure finie $\mu$ sur
$\R^*_+$ telle que \linebreak $\int \la \; d\mu(\la) < \infty$. On y associe la
repr{\'e}sentation unitaire $(U_t)$ de $\R$
sur $L^2(\R^*_+,\mu)$, d{\'e}fini par $(U_t \xi)(\la) = \la^{it}
\xi(\la)$. On suppose que $(U_t)$ est non-p{\'e}riodique.

Connes d{\'e}finit $P =
M_2(L^\infty(\R^*_+,\mu))$ muni de l'{\'e}tat $\vfi$ proportionnel {\`a} la
forme positive
$$\om \begin{pmatrix} f_{11} & f_{12} \\ f_{21} & f_{22} \end{pmatrix}
=
\int f_{11}(\la) \; d\mu(\la) + \int \la f_{22}(\la) \; d\mu(\la) \;
.$$
Prenons un groupe discret infini $G$ et d{\'e}finissons le produit tensoriel infini
$$P_\infty = \bigotimes_{g \in G} (P,\vfi) \; .$$
Alors $G$ agit sur $P_\infty$ par automorphismes de d{\'e}calage des facteurs tensoriels et on
consid{\`e}re le produit crois{\'e} $M = P_\infty \rtimes G$. De cette mani{\`e}re
$M$ est un facteur de type III$_1$. Connes d{\'e}montre que pour $G =
\F_n$ ($n=2,\ldots,+\infty$), $M$ est un facteur plein et l'invariant $\tau(M)$ est la topologie la plus
faible qui rend continue l'application $t \mapsto U_t$. Les facteurs de type III$_1$ de Connes et ceux
de Shlyakhtenko peuvent-ils {\^e}tre isomorphes~?
\end{rema}

\section{Classification des facteurs d'Araki-Woods libres} \label{sec.classification}

\subsection{Le cas presque-p{\'e}riodique}

Supposons d'abord que $(U_t)$ est une repr{\'e}sentation orthogonale \emph{presque-p{\'e}riodique}. Ceci veut dire que l'op{\'e}rateur $A$, qui {\'e}tait d{\'e}fini sur le
complexifi{\'e} $H$ de $H_\R$ par $U_t = A^{it}$, a un spectre purement
ponctuel. Soit $G \subset \R^*_+$
le sous-groupe engendr{\'e} par le spectre ponctuel de
$A$. Shlyakhtenko \cite{shl1} d{\'e}montre que ce sous-groupe \emph{classifie les facteurs d'Araki-Woods libres
  presque-p{\'e}riodiques}.

\begin{theo} \label{theo.classquasiper}
Soit $(U_t)$ une repr{\'e}sentation orthogonale presque-p{\'e}riodique et non-triviale. Soit $G$
le sous-groupe de $\R^*_+$ engendr{\'e} par le spectre ponctuel de
$A$. Alors, \linebreak $(\free{H_\R,U_t},\vfi_U)$ ne d{\'e}pend que de $G$ {\`a} des
isomorphismes qui pr{\'e}servent l'{\'e}tat quasi-libre libre pr{\`e}s.

R{\'e}ciproquement, le groupe $G$ co{\"\i}ncide avec l'invariant $S_{\text{\rm
    discret}}$ du facteur
$\free{H_\R,U_t}$ \cite{Connes-periodic}, qui classifie donc les facteurs d'Araki-Woods libres presque-p{\'e}riodiques et non-triviales.
\end{theo}

En particulier, il d{\'e}coule de ce th{\'e}or{\`e}me et du th{\'e}or{\`e}me \ref{theo.types}
que $(T_\la,\vfi_\la)$ est le seul facteur d'Araki-Woods libre de type
III$_\la$ ($0 < \la < 1$).

Remarquons que le cas o{\`u} $U_t = \id$ pour
tout $t \in \R$ reste ouvert. En effet, d'apr{\`e}s l'isomorphisme
\eqref{eq.trivial}, on sait qu'on obtient le facteur du groupe libre {\`a}
$n$ g{\'e}n{\'e}rateurs~: d{\'e}cider si ces facteurs d{\'e}pendent de $n$ est un des
probl{\`e}mes ouverts en alg{\`e}bres d'op{\'e}rateurs.

\subsection{Le facteur de type II$_\infty$ associ{\'e}}

Au \S\ref{sec.rappels} nous avons vu qu'on associe {\`a} chaque alg{\`e}bre de
von Neumann $M$, une alg{\`e}bre de von Neumann semi-finie $M \rtimes_{(\si_t)}
\R$ o{\`u} $(\si_t)$ est le groupe modulaire d'un {\'e}tat fid{\`e}le sur $M$.
On sait que $M$ est un facteur de type III$_1$ ssi le produit crois{\'e}
$M \rtimes_{(\si_t)} \R$ est un facteur de type II$_\infty$. On
l'appelle le facteur II$_\infty$ associ{\'e} au facteur $M$ de type III$_1$.

Si $M$ est un facteur de type III$_\la$ ($0 < \la < 1$), on peut
prendre un {\'e}tat fid{\`e}le sur $M$ tel que le groupe modulaire
correspondant $(\si_t)$ admette
$\frac{2 \pi}{|\log \la|}$ comme p{\'e}riode. Le groupe modulaire donne donc une action
du cercle $\T$ sur $M$. Le produit crois{\'e} $M \rtimes \T$ est un
facteur de type II$_\infty$ et on a un isomorphisme canonique
$$M \rtimes_{(\si_t)} \R \cong (M \rtimes \T) \ot L^\infty(\T) \; .$$
Le facteur $M \rtimes \T$ de type II$_\infty$ s'appelle {\'e}galement le facteur
II$_\infty$ associ{\'e} {\`a} $M$.

Les facteurs II$_\infty$ associ{\'e}s aux facteurs de type III$_\la$ ($0 <
\la \leq 1$) retiennent une certaine partie de la structure de
$M$. Dans ce paragraphe on s'en sert pour d{\'e}montrer certains r{\'e}sultats
de \emph{non-isomorphisme} entre les facteurs d'Araki-Woods libre.

Dans \cite{shl1} Shlyakhtenko calcule le facteur II$_\infty$ associ{\'e}
au facteur d'Araki-Woods libre $(T_\la,\vfi_\la)$ de type III$_\la$. Il
identifie, gr{\^a}ce aux mod{\`e}les matriciels, le facteur
$T_\la$ au facteur suivant {\'e}tudi{\'e} par R{\u a}dulescu \cite{radulescu-core}.
$$D_\la := (M_2(\C),\om_\la) \ast (L^\infty[-1,1],\mu)
\; ,$$
o{\`u} $\om_\la(e_{ij}) = \delta_{ij} \frac{\la^j}{1 + \la}$ pour $i,j =
0,1$, et $\mu$ est la mesure semi-circulaire. Dans \cite{radulescu-core} R{\u
  a}dulescu d{\'e}montre que le facteur II$_\infty$ associ{\'e} {\`a} $D_\la$ est
isomorphe {\`a} $L(\F_\infty) \ot \B(\ell^2)$. On obtient donc le r{\'e}sultat suivant.

\begin{prop}
Le facteur II$_\infty$ associ{\'e} au facteur d'Araki-Woods libre
$(T_\la,\vfi_\la)$ est isomorphe {\`a} $L(\F_\infty) \ot \B(\ell^2)$.
\end{prop}

Dans \cite{shl3,shl2} une description plus syst{\'e}matique des facteurs
II$_\infty$ associ{\'e}s aux facteurs d'Araki-Woods libres est
donn{\'e}e. L'id{\'e}e est la suivante~: dans l'isomorphisme \eqref{eq.trivial}
nous avons vu que le facteur d'un groupe libre est engendr{\'e} par une
famille libre d'op{\'e}rateurs semi-circulaires. Une telle famille peut
{\^e}tre obtenue par des op{\'e}rateurs de cr{\'e}ation sur un espace de Fock plein.

Shlyakhtenko g{\'e}n{\'e}ralise
ceci et consid{\`e}re dans \cite{shl3} une famille libre
  d'op{\'e}rateurs semi-circulaires \emph{{\`a} coefficients dans une alg{\`e}bre de von Neumann
  $A$}. On retrouve le cas pr{\'e}c{\'e}dent quand $A = \C$. On peut
construire une telle famille {\`a} coefficients dans $A$ en rempla{\c c}ant,
dans la construction de l'espace de Fock plein, les espaces de Hilbert
par des $A$-modules Hilbertiens. Ceci permet d'engendrer le facteur
II$_\infty$ associ{\'e} {\`a} un facteur d'Araki-Woods libre par une famille
libre {\`a} coefficients dans $A=L^\infty(\R)$.

De cette mani{\`e}re Shlyakhtenko d{\'e}montre dans \cite{shl5,shl2} le
r{\'e}sultat suivant.

\begin{theo}
Soit $(H_\R,U_t)$ un multiple fini ou infini de la repr{\'e}sentation r{\'e}guli{\`e}re de $\R$ donn{\'e}e par $(L^2(\R,\R),\la_t)$. Alors, le facteur II$_\infty$
associ{\'e} {\`a} $\free{H_\R,U_t}$ est isomorphe {\`a}  $L(\F_\infty) \ot \B(\ell^2)$.
\end{theo}

Dans \cite{shl2} Shlyakhtenko identifie l'action duale sur le facteur $L(\F_\infty)
\ot \B(\ell^2)$ de type II$_\infty$
associ{\'e} {\`a} $\free{L^2(\R,\R),\la_t}$ avec l'action construite par R{\u
  a}dulescu dans \cite{radulescu-tracescaling}.

\subsection{Des r{\'e}sultats de non-isomorphisme} \label{sec.noniso}

Comme on a vu au \S\ref{sec.rep}, on peut associer {\`a}
chaque mesure sym{\'e}trique $\mu$ sur $\R$, une repr{\'e}sentation
orthogonale $(U_t)$ de $\R$ sur l'espace de Hilbert r{\'e}el $H_\R$ d{\'e}fini
par
$$H_\R = \{\xi \in L^2(\R,\mu) \mid \xi(-x) =
\overline{\xi(x)} \} \quad\text{et}\quad (U_t \xi)(x) = e^{itx} \xi(x)
$$
On notera $\tau(\mu)$ la topologie la plus
faible sur $\R$ qui rend continue l'application $t \recht U_t$ de $\R$
dans $O(H_\R)$ muni de la topologie faible. D'apr{\`e}s le th{\'e}or{\`e}me
\ref{theo.types}, $\tau(\mu)$ est exactement l'invariant $\tau$ du
facteur d'Araki-Woods libre $\free{H_\R,U_t}$. Dans \cite{shl5}
Shlyakhtenko d{\'e}montre qu'il existe une famille non-d{\'e}nombrable de
mesures $\mu$ sans atomes telles que les topologies $\tau(\mu)$ soient
distinctes. On obtient le r{\'e}sultat suivant.

\begin{prop}
Il existe une famille non-d{\'e}nombrable de facteurs d'Araki-Woods libres
mutuellement non-isomorphes et sans {\'e}tats presque-p{\'e}riodiques.
\end{prop}

Les alg{\`e}bres de cette famille peuvent {\^e}tre distingu{\'e}s par l'invariant $\tau$. N{\'e}anmoins, on verra plus tard que l'invariant $\tau$ ne suffit pas pour
distinguer tous les facteurs d'Araki-Woods libres.

Voiculescu a introduit \cite{Voic-entropy2} la notion \emph{d'entropie libre}
$\chi(x_1,\ldots,x_n)$ pour des {\'e}l{\'e}ments auto-adjoints $x_1,\ldots,x_n$
dans une alg{\`e}bre de von Neumann finie $M$ munie d'une trace. Ceci est utilis{\'e}
pour d{\'e}finir la \emph{dimension entropique libre}
$\delta(x_1,\ldots,x_n)$. Une application spectaculaire de l'entropie
libre a {\'e}t{\'e} donn{\'e} par Voiculescu dans \cite{Voic-entropy3} o{\`u} il d{\'e}montre que
les facteurs des groupes libres n'admettent pas de sous-alg{\`e}bre de
Cartan.

Dans \cite{shl5,shl3} Shlyakhtenko utilise la dimension entropique
libre pour d{\'e}montrer que, dans certains cas, le facteur II$_\infty$
associ{\'e} {\`a} un facteur d'Araki-Woods libre ne peut {\^e}tre isomorphe {\`a}
$L(\F_\infty) \ot \B(\ell^2)$.

\begin{theo}
Soit $(U_t)$ une repr{\'e}sentation orthogonale non-p{\'e}riodique de $\R$ sur un espace de
Hilbert r{\'e}el $H_\R$. Supposons que la mesure spectrale de
$\bigoplus_{n \geq 1} U_t^{\ot n}$ est singuli{\`e}re par rapport {\`a} la
mesure de Lebesgue. Alors, le facteur II$_\infty$ associ{\'e} {\`a}
$\free{H_\R,U_t}$ n'est pas isomorphe {\`a} $L(\F_\infty) \ot \B(\ell^2)$.

En particulier, $\free{H_\R,U_t}$ n'est pas isomorphe {\`a} $\free{L^2(\R,\R),\la_t}$, o{\`u} $(\la_t)$ est la repr{\'e}\-sentation r{\'e}guli{\`e}re de $\R$.
\end{theo}

La condition du th{\'e}or{\`e}me pr{\'e}c{\'e}dent est satisfaite si la topologie
la plus faible qui rend continue l'application $t \mapsto U_t$ est
strictement plus faible que la topologie usuelle de $\R$.

Dans \cite{shl5} Shlyakhtenko construit une mesure $\mu$ sur $\R$
telle que toutes les mesures $\mu \ast \cdots \ast \mu$ sont
singuli{\`e}res par rapport {\`a} la mesure de Lebesgue, mais n{\'e}anmoins
$\tau(\mu)$ est la topologie usuelle de $\R$. Le th{\'e}or{\`e}me pr{\'e}c{\'e}dent admet
donc le corollaire suivant.

\begin{coro} \label{coro.indistinguable}
Il existe des facteurs d'Araki-Woods libres non-isomorphes ayant le m{\^e}me invariant $\tau$.
\end{coro}

Comme l'invariant $\tau$ ne distingue pas tous
les facteurs d'Araki-Woods libres, Shlyakhtenko propose dans \cite{shl5}
un \emph{nouvel invariant} $\cS$ pour les facteurs pleins de type III.
Introduisons quelques notations.
Si $\mu$ est une mesure sur $\R$, notons
$\cC_\mu$ l'ensemble de toutes les mesures qui sont
absolument continues par rapport {\`a} la mesure $\mu$. Dans le cas o{\`u}
$\mu$ est la mesure spectrale d'un op{\'e}rateur auto-adjoint et
strictement positif $A$, on pose $\cC_A := \cC_\mu$
Ceci permet de d{\'e}finir
$$\cS(M) := \bigcap_{\vfi \; \text{{\'e}tat fid{\`e}le sur} \; M}
\cC_{\bigoplus_n \Delta_\vfi^{\ot n}} \; ,$$
o{\`u} $\Delta_\vfi$ est l'op{\'e}rateur modulaire de l'{\'e}tat $\vfi$. On
remarque que les mesures dans $\cS(M)$ sont support{\'e}es par $\R^*_+$ et
que la mesure de Dirac $\delta_1$ est toujours dans $\cS(M)$.

Shlyakhtenko d{\'e}montre dans \cite{shl5} que cet invariant $\cS$ distingue certains
facteurs d'Araki-Woods libres (non-isomorphes) qui ont le m{\^e}me invariant $\tau$.

\subsection{Le facteur d'Araki-Woods libre d{\'e}pend-il de la multiplicit{\'e}~?}

A chaque mesure sym{\'e}trique $\mu$ sur $\R$ est associ{\'e}e une repr{\'e}sentation orthogonale (voir \S\ref{sec.noniso}). Notons $\Gamma(\mu,n)$ o{\`u} $n \in
\{1,\ldots,+\infty\}$, le facteur d'Araki-Woods libre associ{\'e} {\`a} la somme directe de $n$ copies de cette repr{\'e}sentation.

Il d{\'e}coule du th{\'e}or{\`e}me de classification \ref{theo.classquasiper} que,
dans le cas o{\`u} $\mu$ est une mesure atomique non-concentr{\'e}e en
$\{0\}$, le facteur $\Gamma(\mu,n)$ ne d{\'e}pend pas de $n$. D'apr{\`e}s le
th{\'e}or{\`e}me \ref{theo.types}, l'invariant $\tau$ d'un facteur
$\Gamma(\mu,n)$ quelconque ne d{\'e}pend pas de $n$. N{\'e}anmoins,
Shlyakhtenko d{\'e}montre dans \cite{shl6} un r{\'e}sultat tr{\`e}s surprenant.

\begin{theo} \label{theo.dependmult}
Soit $\la$ la mesure de Lebesgue sur $\R$ et $\delta_0$ la mesure de
Dirac en $0$. Alors, $\Gamma(\la+\delta_0,1)$ et
$\Gamma(\la+\delta_0,2)$ ne sont pas isomorphes.
\end{theo}

Dans sa preuve Shlyakhtenko utilise la notion \emph{d'alg{\`e}bre de von
  Neumann solide} due {\`a} Ozawa \cite{ozawa}~: une alg{\`e}bre de von Neumann
  est dite solide si le commutant relatif de n'importe quelle
  sous-alg{\`e}bre diffuse et unif{\`e}re est injectif. Rappelons qu'une
  alg{\`e}bre de von Neumann est dite diffuse si elle n'admet pas de
  projecteurs minimaux. Une alg{\`e}bre de von Neumann solide est
  n{\'e}cessairement finie.

Ozawa d{\'e}montre dans \cite{ozawa} que l'alg{\`e}bre de von Neumann $L(G)$
d'un groupe discret hyperbolique $G$ (voir \cite{ghys-etal}) est solide. En
particulier, les facteurs des groupes libres sont solides.

Notons $N_n$ le facteur II$_\infty$ associ{\'e} {\`a} $\Gamma(\la+\delta_0,n)$. Shlyakhtenko d{\'e}montre que $N_1 \cong L(\F_\infty) \ot
\B(\ell^2)$. Ceci implique que $pN_1p$ est une alg{\`e}bre de von Neumann
solide pour tout projecteur fini $p \in N_1$. Par contre, il
construit {\'e}galement un projecteur fini $q \in N_2$ tel que $qN_2q$ ne
soit pas solide.

Remarquons qu'il d{\'e}coule des r{\'e}sultats de \cite{shl5} que l'invariant $\cS$ ne distingue pas \linebreak $\Gamma(\la+\delta_0,1)$ et
$\Gamma(\la+\delta_0,2)$.

Soit $(U_t)$ une repr{\'e}sentation orthogonale qui contient une repr{\'e}sentation p{\'e}riodique
non-triviale.
D'apr{\`e}s le th{\'e}or{\`e}me \ref{theo.absorption} (et la proposition
\ref{prop.sommedirecte}), on sait que le facteur d'Araki-Woods libre
associ{\'e} absorbe
$(L(\F_\infty),\trace)$~:
$$(\free{H_\R,U_t},\vfi_U) \cong (\free{H_\R,U_t},\vfi_U) \ast
(L(\F_\infty),\trace) \; .$$
Le deuxi{\`e}me r{\'e}sultat surprenant de \cite{shl6} est qu'il
existe des facteurs d'Araki-Woods libres qui n'absorbent pas
$(L(\F_\infty),\trace)$.

\begin{theo}
Soit $\la$ la mesure de Lebesgue sur $\R$. Alors,
$$\Gamma(\la,1) \not\cong (\Gamma(\la,1),\vfi_{\la,1}) \ast
(L(\F_\infty),\trace) \; .$$
\end{theo}

\section{Appendice : sur le lemme des $\mathbf{14}\mbox{\large $\eps$}$} \label{sec.appendice}

Dans \ref{lemm.technique} nous d{\'e}montrons une g{\'e}n{\'e}ralisation du lemme
technique 4.1 de \cite{Vouter}, qui {\'e}tait {\`a} son tour une
g{\'e}n{\'e}ralisation du
lemme des $14\eps$ d{\^u} {\`a} Murray \& von Neumann \cite{murr-vn} (voir
\cite{Bar} pour une version adapt{\'e}e aux produits libres de type III). On exploite
le fait que le produit libre $G_1 \ast G_2$ de deux groupes non-triviaux $G_1,G_2$ est
tr{\`e}s non-moyennable\footnote{Plus pr{\'e}cis{\'e}ment $G_1 \ast G_2$ n'est pas
int{\'e}rieurement moyennable. On construit {\'e}galement une d{\'e}composition paradoxale
explicite de $G$.}, sauf si $G_1 \cong G_2 \cong \Z / 2\Z$. Ceci explique
le lemme \ref{lemm.technique}~: il nous faut un {\'e}l{\'e}ment non-trivial
dans $N_1$ et deux {\'e}l{\'e}ments non-triviaux dans $N_2$.

Le lemme
\ref{lemm.technique} permet de calculer l'invariant $\tau$ d'un
certain nombre de produits libres, voir
proposition \ref{prop.invarianttau}.

\begin{lemm} \label{lemm.technique}
Soit $N_i$ une alg{\`e}bre de von Neumann munie d'un {\'e}tat fid{\`e}le
$\om_i$, $(i=1,2)$. Posons $(N,\om) = (N_1,\om_1) \ast
(N_2,\om_2)$. Soient $a \in N_1$ et
$b,c \in N_2$. Supposons que les {\'e}l{\'e}ments $a,b$ et $c$ appartiennent
au domaine de $\si^\om_{i/2}$, o{\`u} $(\si^\om_t)$ est le groupe
modulaire de
l'{\'e}tat $\om$. Soit $\al_i$ un automorphisme de $N_i$ qui satisfait $\om_i \al_i = \om_i$, $(i=1,2)$. Notons $\al = \al_1 \ast \al_2$. Alors, pour
tout $x \in N$,
\begin{align*}
\|x - & \om(x)1 \|_2 \leq \cE(a,b,c)\, \max \bigl\{ \|xa - \al(a)x\|_2 , \|xb - \al(b)x\|_2,\|xc - \al(c)x\|_2 \bigr\} \\ & \hspace{12cm} + \cF(a,b,c) \, \|x\|_2 \\
\text{o{\`u}}\;\; &\cE(a,b,c) = 6\|a\|^3 + 4 \|b\|^3 + 4 \|c\|^3 \; , \\
& \cF(a,b,c) = 3 \cC(a) + 2 \cC(b) + 2 \cC(c) + 12 |\om(cb^*)| \, \|cb^*\| \; ,\\
& \cC(a) = 2\|a\|^3 \, \|\si^\om_{i/2}(a) - a \| + 2\|a\|^2 \, \|a^*a-1\| + 3(1+\|a\|^2)\,\|aa^*-1\| + 6 |\om(a)| \, \|a\| \; .
\end{align*}
\end{lemm}

\begin{proof}
Repr{\'e}sentons $N_i$ sur l'espace de Hilbert $H_i$ de la repr{\'e}sentation
 GNS de $\om_i$ et soit $\xi_i$ le vecteur cyclique associ{\'e}. Posons
$(H,\xi) = (H_1,\xi_1) \ast (H_2,\xi_2)$. On rappelle \cite{Voic-livre} que
$$H= \C \xi \oplus (\Honebol \ot H(2,l)) \oplus (\Htwobol \ot H(1,l))
\; ,$$ o{\`u} $\Hbol_i = H_i \ominus \C \xi_i$,
\begin{align*}
& H(2,l) = \C \xi \oplus \Htwobol \oplus (\Htwobol \ot \Honebol) \oplus (\Htwobol \ot \Honebol \ot \Htwobol) \oplus \cdots  \; , \\ & H(1,l) = \C \xi
\oplus \Honebol \oplus (\Honebol \ot \Htwobol) \oplus (\Honebol \ot \Htwobol \ot \Honebol) \oplus \cdots \; .
\end{align*}
Pour $\zeta \in H$ et $y \in N$, on d{\'e}finit l'action {\`a} droite de $y$ sur $\zeta$ par $\zeta \cdot y := Jy^*J \zeta$ o{\`u} $J$ est la conjugaison
modulaire de l'{\'e}tat $\om$.

Choisissons $x \in N$ et d{\'e}finissons $\eta = x \xi$. On {\'e}crit $\eta = \om(x) \xi + \mu + \ga$ avec $\mu \in \Honebol \ot H(2,l)$ et $\ga \in \Htwobol
\ot H(1,l)$. Posons alors $\xbol = x - \om(x)1$, $\eta_0 = \mu + \ga$, $\etatil = \al(a^*) \cdot \eta \cdot a$, $\gatil = \al(a^*) \cdot \ga \cdot a$
et $\zetatil = \eta_0 - \ga - \gatil$. Bien {\'e}videmment
\begin{align*}
\|\mu\|^2 + \|\ga\|^2 &= \|\zetatil + \ga + \gatil\|^2 \geq \|\zetatil\|^2 + \|\ga\|^2 + \|\gatil\|^2 - 2 |\langle \zetatil,\ga\rangle| - 2 |\langle
\zetatil,\gatil \rangle| - 2 |\langle \ga,\gatil \rangle|
\\ & \geq 2\|\ga\|^2 - \bigl| \|\ga\|^2 - \|\gatil\|^2\bigr| - 2 |\langle
\zetatil,\ga\rangle| - 2 |\langle \zetatil,\gatil \rangle| - 2 |\langle \ga,\gatil \rangle| \; .
\end{align*}
Exactement de la m{\^e}me mani{\`e}re que dans la d{\'e}monstration du lemme 4.1 de \cite{Vouter}, on sait estimer tous les termes n{\'e}gatifs. On conclut que
\begin{equation}\label{eq.firstfinal} \|\ga\|^2 \leq \|\mu\|^2 + 2 \|a\|^3 \; \|xa - \al(a)x \|_2 \; \|\xbol\|_2 + \cC(a) \; \|x\|_2 \; \|\xbol\|_2
\; .
\end{equation}
On obtient une estimation analogue {\`a} l'aide des {\'e}l{\'e}ments $b$ et $c$.
En effet on pose $\etapr = \al(b^*)
\cdot \eta \cdot b$, $\etadpr = \al(c^*) \cdot \eta \cdot
c$, $\mupr = \al(b^*) \cdot \mu \cdot b$ et $\mudpr = \al(c^*) \cdot \mu \cdot
c$. On d{\'e}finit $\zetapr = \eta_0 - \mu - \mupr -
\mudpr$. On trouve que
\begin{align*}
\|\mu\|^2 + \|\ga\|^2 \geq
& 3 \|\mu\|^2 - \bigl| \|\mu\|^2 - \|\mupr\|^2 \bigr| -
\bigl| \|\mu\|^2 - \|\mudpr\|^2 \bigr|
- 2 |\langle \zetapr,\mu \rangle| - 2 |\langle \zetapr,\mupr \rangle|
\\ &- 2 |\langle \zetapr,\mudpr \rangle|  - 2 |\langle \mu,\mupr \rangle|
- 2 |\langle \mu,\mudpr \rangle|  - 2 |\langle \mupr,\mudpr \rangle|
\; .
\end{align*}
On estime de nouveau tous les termes n{\'e}gatifs et on obtient
\begin{equation} \label{eq.secondfinal}
\begin{split}
2 \|\mu\|^2 \leq \|\ga\|^2 & + 2 \|b\|^3 \; \|xb - \al(b)x \|_2 \; \|\xbol\|_2 + 2 \|c\|^3 \; \|xc - \al(c)x \|_2 \; \|\xbol\|_2 \\ &+\bigl( \cC(b) +
\cC(c) + 6 \|cb^*\| \; |\om(cb^*)| \bigr) \; \|x\|_2 \; \|\xbol\|_2 \; .
\end{split}
\end{equation}
Comme $\|\mu\|^2 + \|\ga\|^2 = \|\xbol\|_2^2$, une combinaisons des in{\'e}galit{\'e}s \eqref{eq.firstfinal} et \eqref{eq.secondfinal} donne l'in{\'e}galit{\'e} du
lemme.
\end{proof}

A priori l'invariant $\tau$ d'un facteur plein est difficile {\`a}
calculer car on munit $\R$ de la topologie induite d'une topologie quotient.
L'int{\'e}r{\^e}t de la proposition
suivante est de donner une formule pour l'invariant $\tau$ en termes
du groupe modulaire d'un seul {\'e}tat fid{\`e}le, sans qu'il faille conna{\^\i}tre $\Out
N$. De la m{\^e}me mani{\`e}re que Shlyakhtenko d{\'e}duit du lemme des $14\eps$ de Barnett son corollaire 8.4 dans \cite{shl3}, nous d{\'e}duisons du lemme
\ref{lemm.technique} le r{\'e}sultat suivant, utilis{\'e} dans le \S\ref{sec.circulaire}.

\begin{prop} \label{prop.invarianttau}
Soit $N_i$ des alg{\`e}bres de von Neumann munies d'un {\'e}tat fid{\`e}le $\om_i$, $(i=1,2)$. Soit $(N,\om) = (N_1,\om_1) \ast (N_2,\om_2)$. On suppose que
$N_1$ contient une suite $(a_n)$ d'{\'e}l{\'e}ments qui sont analytique par rapport {\`a} l'{\'e}tat $\om$ et satisfont
\begin{equation}\label{eq.suite}
\|\si^\om_{i/2}(a_n) - a_n\| \recht 0 \; , \quad \|a_n^* a_n - 1\|, \|a_n a_n^* - 1\| \recht 0 \quad\text{et}\quad \om(a_n) \recht 0 \; .
\end{equation}
On suppose que $N_2$ contient des suites $(b_n),(c_n)$ qui satisfont les m{\^e}mes conditions que $(a_n)$ ainsi que la condition $\om(c_n b_n^*) \recht
0$. Alors,
\begin{itemize}
\item[a)] $N$ est un facteur plein.
\item[b)] Notons $\Aut(N_i,\om_i)$ le groupe d'automorphismes de $N_i$
  pr{\'e}servant l'{\'e}tat $\om_i$ et $\pi : \Aut(N) \recht \Out(N)$
  l'application quotient. Alors,
l'homomorphisme
\begin{equation}\label{eq.homom}
\Aut(N_1,\om_1) \times \Aut(N_2,\om_2) \recht
  \Out(N) : (\al_1,\al_2) \mapsto \pi(\al_1 \ast \al_2)
\end{equation}
est un
  hom{\'e}omorphisme {\`a} image ferm{\'e}.
\item[c)] L'invariant $\tau(N)$ est la topologie la plus faible qui rend continues les applications $t \mapsto \si^{\om_i}_t$
de $\R$ dans $\Aut N_i$ ($i=1,2$).
\end{itemize}
\end{prop}

\begin{proof}
Soit $(x_k)$ une suite d'unitaires dans $N$ et $\al_k \in
\Aut(N_1,\om_1)$, $\be_k \in \Aut(N_2,\om_2)$ des
suites d'automorphismes. Supposons que $\Ad(x_k^*) \circ (\al_k \ast
\be_k) \recht \id$ dans $\Aut(N)$. Il suffit de d{\'e}montrer que $\|x_k -
\om(x_k) 1 \|_2 \recht 0$. En effet, prenant $\al_k = \id$ et $\be_k =
\id$ pour tout $k$, on aura d{\'e}montr{\'e} que chaque suite centrale est
triviale et que donc, $N$ est un facteur plein. On aura {\'e}galement
d{\'e}montr{\'e} que l'homomorphisme \eqref{eq.homom} est un
hom{\'e}omorphisme. Finalement c) r{\'e}sulte de b).

Choisissons $\eps > 0$. Prenons $n$ tel que $\cF(a_n,b_n,c_n) <
\eps/2$. Si $k \recht \infty$, on a
$$\|x_k a_n - (\al_k \ast \be_k)(a_n) x_k \|_2 = \|(\id - \Ad(x_k^*) \circ (\al_k \ast
\be_k))(a_n) \|_2 \recht 0$$
et on a le m{\^e}me r{\'e}sultat rempla{\c c}ant $(a_n)$ par $(b_n)$ ou $(c_n)$.
Il d{\'e}coule du lemme \ref{lemm.technique} qu'il existe $n_0$ tel que
pour tout $n \geq n_0$ on a $\|x_n - \om(x_n)1 \|_2 < \eps$.
\end{proof}

On peut appliquer la proposition \ref{prop.invarianttau} {\`a} un produit
libre $(N_1,\om_1) \ast \bigl( (N_2,\om_2) \ast (N_3,\om_3) \bigr)$ si
chacune des alg{\`e}bres $(N_i,\om_i)$ contient une suite $(a_n)$ qui
satisfait les conditions \eqref{eq.suite}. En effet, comme $b_n$ et
$c_n$ sont $^*$-libres dans ce cas-l{\`a}, on a automatiquement que $\om(c_n b_n^*) \recht
0$. En particulier, il d{\'e}coule de la proposition
\ref{prop.sommedirecte} qu'on peut appliquer la proposition
\ref{prop.invarianttau} {\`a} un facteur d'Araki-Woods libre associ{\'e} {\`a} une
repr{\'e}sentation orthogonale de $\R$ qui est une somme directe de trois
repr{\'e}sentations, pourvu qu'on d{\'e}montre que chaque facteur d'Araki-Woods libre contient une
suite $(a_n)$ qui satisfait les conditions \eqref{eq.suite}.

\begin{lemm} \label{lemm.approx}
Soit $(U_t)$ une repr{\'e}sentation orthogonale de $\R$ sur l'espace de Hilbert r{\'e}el $H_\R$. Soit $\free{H_\R,U_t}$ l'alg{\`e}bre de von
Neumann associ{\'e}e {\`a} l'{\'e}tat quasi-libre libre $\vfi_U$. Soit $(\si_t)$ le groupe modulaire de l'{\'e}tat $\vfi_U$. Alors il existe une suite d'unitaires
$(u_n)$ dans $\free{H_\R,U_t}$ qui sont analytiques par rapport {\`a} $(\si_t)$ et satisfont
$$\|\si_z(u_n) - u_n \| \recht 0 \;\;\text{uniform{\'e}ment sur des compacts de $\C$, et}\quad \vfi_U(u_n) \recht 0 \; .$$
\end{lemm}

\begin{proof}
Si $U_t = \id$ pour tout $t \in \R$, le lemme est trivial. On suppose donc que $(U_t)$ est non-trivial. Soit $A$ l'op{\'e}rateur auto-adjoint strictement
positif sur le complexifi{\'e} $H$ de $H_\R$ tel que $U_t = A^{it}$. Soit $J$ l'anti-unitaire canonique de $H$ et $T = JA^{-1/2}$ l'involution sur $H$
associ{\'e}e {\`a} $(U_t)$. Comme $A \neq 1$ et $JAJ = A^{-1}$, on peut prendre $\lambda > 1$ dans le spectre de $A$. Notons $\chi_n$ la fonction
indicatrice de l'intervalle $[\lambda-\frac{1}{n},\lambda+\frac{1}{n}]$ et prenons des vecteurs unit{\'e} $\xi_n$ dans l'image de $\chi_n(A)$. Comme $JAJ
= A^{-1}$, les vecteurs $\xi_n$ et $J\xi_n$ seront orthonormaux pour $n$ suffisamment grand. On d{\'e}finit les {\'e}l{\'e}ments $x_n \in \free{H_\R,U_t}$ par
$x_n = \ell(\xi_n) + \ell(T\xi_n)^*$. Il est clair que $x_n$ est analytique par rapport {\`a} $(\si_t)$ et que $\| \si^\om_z(x_n) - \lambda^{iz} x_n \|
\recht 0$ uniform{\'e}ment sur des compacts de $\C$.

D{\'e}finissons l'op{\'e}rateur $y_n = \ell(\xi_n) + \frac{1}{\sqrt{\lambda}} \ell(J \xi_n)^*$ dans $\B(\cF(H))$. Alors $\|x_n - y_n \| \recht 0$ et d'apr{\`e}s
le th{\'e}or{\`e}me \ref{theo.Tla}, l'op{\'e}rateur $y_n^* y_n$ a une distribution sans atomes par rapport {\`a} l'{\'e}tat vectoriel du vide qui ne d{\'e}pend pas de $n$.
Il existe donc une function continue $g : \R \recht \R$ telle que $\langle \exp(i g(y_n^* y_n)) \Omega,\Omega \rangle=0$ pour tout $n$.

Choisissons $\eps > 0$ et $K \subset \C$ compact. Comme $g$ peut {\^e}tre approxim{\'e}e par des polyn{\^o}mes uniform{\'e}ment sur le spectre de $y_n^*y_n$, on peut
prendre un polyn{\^o}me $P$ avec des coefficients r{\'e}els tel que $|\langle \exp(i p(y_n^* y_n)) \Omega,\Omega \rangle| < \eps/2$ pour tout $n$. On sait
que $\|\si^\om_z(x_n^* x_n) - x_n^* x_n\| \recht 0$ uniform{\'e}ment sur des compacts de $\C$ et que $\|x_n - y_n\| \recht 0$. Pour $n$ suffisamment
grand $u := \exp(iP(x_n^* x_n))$ est alors un unitaire dans $\free{H_\R,U_t}$ qui satisfait $\|\si^\om_z(u) - u \| < \eps$ pour tout $z \in K$ et
$|\vfi_U(u)| < \eps$.
\end{proof}

\end{document}